\documentclass[a4paper,10pt]{article}
\usepackage{cstafa}

\title{Models of Non-Well-Founded Sets via an Indexed Final Coalgebra Theorem}

\author{Benno van den Berg \and Federico De Marchi}

\begin{document}

\maketitle

\begin{abstract}
  \noindent The paper uses the formalism of indexed categories to
  recover the proof of a standard final coalgebra theorem, thus
  showing existence of final coalgebras for a special class of
  functors on categories with finite limits and colimits. As an
  instance of this result, we build the final coalgebra for the
  powerclass functor, in the context of a Heyting pretopos with a
  class of small maps. This is then proved to provide a model for
  various non-well-founded set theories, depending on the chosen
  axiomatisation for the class of small maps.
\end{abstract}

\section{Introduction}

\begin{quote}
  {\it The explicit use of bisimulation for set theory goes back to
  the work on non-wellfounded sets by Aczel(1988). It would be of
  interest to construct sheaf models for the theory of non-wellfounded
  sets from our axioms for small maps.}
  \flushright{ -- Joyal and Moerdijk, 1995}
\end{quote}

Since its first appearance in the book by Joyal and Moerdijk
\cite{joyalmoerdijk95}, algebraic set theory has always claimed the
virtue of being able to describe, in a single framework, various
different set theories. In fact, the correspondence between axiom
systems for a class of small maps and formal set theories has been put
to work first in the aforementioned book, and then in the work by
Awodey et al. \cite{awodeybutzsimpsonstreicher04}, thus modelling such
theories as {\bf CZF}, {\bf IZF}, {\bf BIST}, {\bf CST} and so
on. However, despite the suggestion in \cite{joyalmoerdijk95}, it
appears that up until now no one ever tried to put small maps to use
in order to model a set theory which includes the Anti-Foundation
Axiom {\bf AFA}.

This papers provides a first step in this direction. In particular, we
build a categorical model of the weak constructive theory {\bf
CZF$_0$} of (possibly) non-well-founded sets, studied by Aczel and
Rathjen in \cite{aczelrathjen01} . Classically, the universe of
non-well-founded sets is known to be the final coalgebra of the
powerclass functor \cite{aczel88}. Therefore, it should come as no
surprise that we can build such a model from the final coalgebra for
the functor \spower\ determined by a class of small maps.

Perhaps more surprising is the fact that such a coalgebra {\em always}
exists. We prove this by means of a final coalgebra theorem, for a
certain class of functors on a finitely complete and cocomplete
category. The intuition that guided us along the argument is a
standard proof of a final coalgebra theorem by Aczel \cite{aczel88}
for set-based functors on the category of classes that preserve
inclusions and weak pullbacks. Given one such functor, he first
considers the coproduct of all small coalgebras, and shows that this
is a weakly terminal coalgebra. Then, he quotients by the largest
bisimulation on it, to obtain a final coalgebra. The argument works
more generally for any functor of which we know that there is a
generating family of coalgebras, for in that case we can take the
coproduct of that family, and perform the construction as above. The
condition of a functor being set-based assures that we are in such a
situation.

Our argument is a recasting of the given one in the internal language
of a category. Unfortunately, the technicalities that arise when
externalising an argument which is given in the internal language can
be off-putting, at times. For instance, the externalisation of
internal colimits forces us to work in the context of indexed
categories and indexed functors. Within this context, we say that an
indexed functor (which turns pullbacks into weak pullbacks) is
small-based when there is a ``generating family'' of coalgebras. For
such functors we prove an indexed final coalgebra theorem. We then
apply our machinery to the case of a Heyting pretopos with a class of
small maps, to show that the functor \spower\ is small-based and
therefore has a final coalgebra. As a byproduct, we are able to build
the M-type for any small map $f$ (i.e.\ the final coalgebra for the
polynomial functor $P_f$ associated to $f$).

For sake of clarity, we have tried to collect as much indexed
category theory as we could in a separate section. This forms the
content of \refsec{inserters}, and we advise the uninterested reader
to skip all the details of the proofs therein. This should not affect
readability of \refsec{finalcoalgebra}, where we prove our final
coalgebra results. Finally, in \refsec{cst} we prove that the final
\spower-coalgebra is a model of the theory {\bf CZF$_0$}+{\bf AFA}.

Our choice to focus on a weak set theory such as {\bf CZF$_0$} is
deliberate, since stronger theories can be modelled simply by adding
extra requirements for the class of small maps. For example, we can
model the theory {\bf CST} of Myhill \cite{myhill75} (plus {\bf AFA}),
by adding the Exponentiation Axiom, or {\bf IZF$^-$}+{\bf AFA} by
adding the Powerset, Separation and Collection axioms from
\cite[p.~65]{joyalmoerdijk95}. And we can force the theory to be
classical by working in a boolean pretopos. This gives a model of {\bf
ZF$^-$}+{\bf AFA}, the theory presented in Aczel's book
\cite{aczel88}, apart from the Axiom of Choice. Finally, by adding
appropriate axioms for the class of small maps, we build a model of
the theory {\bf CZF$^-$}+{\bf AFA}, which was extensively studied by
M.~Rathjen in \cite{rathjen03,rathjen04}.

As a final remark, we would like to point out that the present results
fit in the general picture described by the two present authors in
\cite{demarchivandenberg04}. (Incidentally, we expect that, together
with the results on sheaves therein, they should yield an answer to
the question by Joyal and Moerdijk which we quoted in opening this
introduction.) There, we suggested that the established connection
between Martin-L\"of type theory, constructive set theory and the
theory of $\Pi W$-pretoposes had an analogous version in the case of
non-well-founded structures. While trying to make the correspondence
between the categorical and the set theoretical sides of the picture
precise, it turned out that the M-types in $\Pi M$-pretoposes are not
necessary, in order to obtain a model of some non-well-founded set
theory. This phenomenon resembles the situation in \cite{lindstrom89},
where Lindstr\"om built a model of {\bf CZF$^-$}+{\bf AFA} out of a
Martin-L\"of type theory with one universe, without making any use of
M-types.

\section{Generating objects in indexed categories}
\label{sec:inserters}

As we mentioned before, our aim is to prove a final coalgebra theorem
for a special class of functors on finitely complete and cocomplete
categories. The proof of such results will be carried out by repeating
in the internal language of such a category \ct{C} a classical
set-theoretic argument. This forces us to consider \ct{C} as an
indexed category, via its canonical indexing \indct{C}, whose fibre
over an object $X$ is the slice category $\ct{C}/X$. We shall then
focus on endofunctors on \ct{C} which are components over $1$ of
indexed endofunctors on \indct{C}. For such functors, we shall prove
the existence of an indexed final coalgebra, under suitable
assumptions. The component over $1$ of this indexed final coalgebra
will be the final coalgebra of the original \ct{C}-endofunctor.

Although in \refsec{fincoals} we will apply our results only in a very
specific setting, it turns out that all the basic machinery needed for
the proofs can be stated in a more general context. This section
collects as much of the indexed category theoretic material as
possible, hoping to leave the other sections easier to follow for a
less experienced reader.

So, for this section, \base{S}\ will be a cartesian category, which we
use as a base for indexing. Our notations for indexed categories and
functors follow those of \cite[Chapters B1 and B2]{johnstone02}, to
which we refer the reader for all the relevant definitions.

We will mostly be concerned with \emph{\base{S}-cocomplete
categories}, i.e.\ \base{S}-indexed categories in which each fibre is
finitely cocomplete, finite colimits are preserved by reindexing
functors, and these functors have left adjoints satisfying the
Beck-Chevalley condition. Under these assumptions it immediately
follows that:
\begin{lemm}{indterm}
  If the fibre $\ct{C}=\fibre{C}{1}$ of an \base{S}-cocomplete
  \base{S}-indexed category \indct{C} has a terminal object $T$, then
  this is an indexed terminal object, i.e.\ $X^*T$ is terminal in
  \fibre{C}{X} for all $X$ in \base{S}.
\end{lemm}

The first step, in the set-theoretic argument to build the final
coalgebra, is to identify a ``generating family'' of coalgebras, in
the sense that any other coalgebra is the colimit of all coalgebras in
that family which map to it. If we want to express this in the
internal language, we need to introduce the concept of internal
colimits in indexed categories. To this end, we first recall that an
{\em internal category} \inct{K} in \base{S} consists of a diagram
\diag{
  K_1\ar@<.8ex>[r]^{d_1} \ar@<-.8ex>[r]_-{d_0} & K_0,
}
where $d_1$ is the \emph{domain} map, $d_0$ is the \emph{codomain} one
and they have a common left inverse $i$, satisfying the usual
conditions. There is also a notion of internal functor between
internal categories, and this gives rise to the category of internal
categories in \base{S}\ (see \cite[Section B2.3]{johnstone02} for the
details).

An \emph{internal diagram $L$ of shape \inct{K}} in an
\base{S}-indexed category \indct{C} consists of an internal
\base{S}-category \inct{K}, an object $L$ in \fibre{C}{K_0}, and a map
$d_1^*L\rTo d_0^*L$ in \fibre{C}{K_1} which interacts properly with
the categorical structure of \inct{K}. Moreover, one can consider the
notion of \emph{morphism of internal diagrams}, and these data define
the category $\indct{C}^\inct{K}$ of \emph{internal diagrams of shape
\inct{K} in \indct{C}}.

An indexed functor \func{F}{\indct{C}}{\indct{D}} induces an ordinary
functor
\func{F^{\inct{K}}}{\indct{C}^{\inct{K}}}{\indct{D}^{\inct{K}}}
between the corresponding categories of internal diagrams of shape
\inct{K}. Dually, given an internal functor
\func{F}{\inct{K}}{\inct{J}}, this (contravariantly) determines by
reindexing of \indct{C} an ordinary functor on the corresponding
categories of internal diagrams:
\func{F^*}{\indct{C}^{\inct{J}}}{\indct{C}^{\inct{K}}}. We say that
\indct{C} has \emph{internal left Kan extensions} if these reindexing
functors have left adjoints, denoted by \Lan{F}{}. In the particular
case where $\inct{J}=1$, the trivial internal category with one
object, we write \func{\inct{K}^*}{\indct{C}}{\indct{C}^{\inct{K}}}
for the functor, and $\colim_{\inct{K}}$ for its left adjoint
\Lan{\inct{K}}{}, and we call $\colim_{\inct{K}}L$ the {\em internal
colimit} of $L$.

\begin{defi}{prescolim}
  Suppose \indct{C} and \indct{D}\ are \base{S}-indexed categories
  with internal colimits of shape \inct{K}. Then, we say that an
  \base{S}-indexed functor \func{F}{\indct{C}}{\indct{D}} {\em
  preserves colimits} if the canonical natural transformation filling
  the square
  \diag{ \indct{C}^{\inct{K}} \ar[d]_{\colim_{\inct{K}}}
  \ar[r]^-{F^{\inct{K}}} & \indct{D}^{\inct{K}}
  \ar[d]^{\colim_{\inct{K}}} \\ \indct{C} \ar[r]_-F & \indct{D}}
  is an isomorphism.
\end{defi}

It follows at once from Proposition B2.3.20 in \cite{johnstone02}
that:

\begin{prop}{indexedcolimits}
  If \indct{C} is an \base{S}-cocomplete \base{S}-indexed category,
  then it has colimits of internal diagrams and left Kan extensions
  along internal functors in \base{S}. Moreover, if an indexed functor
  \func{F}{\indct{C}}{\indct{D}} between \base{S}-cocomplete
  categories preserves \base{S}-indexed colimits, then it also
  preserves internal colimits.
\end{prop}

When forming the internal diagram of those coalgebras that map into a
given one, say $(A,\alpha)$, we need to select out of an object of
maps to $A$ those which are coalgebra morphisms. In order to consider
such objects of arrows in the internal language, we need to introduce
the following concept:

\begin{defi}{exponentiable}
 An object $E$ in the fibre \fibre{C}{U} of an \base{S}-indexed
 category \indct{C} is called \emph{exponentiable} if for any object
 $A$ in any fibre \fibre{C}{I} there is an \emph{exponential} $A^E$
 fitting in a span
 \diaglab{exponential}{U& \ar[l]_-s A^E \ar[r]^-t & I}
 in \base{S} and a {\em generic arrow} \func{\varepsilon}{s^*E}{t^*A} in
 \fibre{C}{A^E}, with the following universal property: for any other
 span in \base{S}
 \diag{U & \ar[l]_-x J \ar[r]^-y &I}
 and any arrow \func{\psi}{x^*E}{y^*A} in \fibre{C}{J}, there is a
 unique arrow \func{\chi}{J}{A^E} in \base{S} such that $s\xi=x$,
 $t\chi=y$ and $\chi^*\varepsilon\cong\psi$ (via the canonical
 isomorphisms arising from the two previous equalities).
\end{defi}

\begin{rema}{expnsstable}
  It follows from the definition, via a standard diagram chasing, that
  the reindexing along an arrow \func{f}{V}{U} in \base{S} of an
  exponentiable object $E$ in \fibre{C}{U} is again exponentiable.
\end{rema}

\begin{rema}{expncanindexing}
  We advise the reader to check that, in case \ct{C} is a cartesian
  category and \indct{C} is its canonical indexing over itself, the
  notion of exponentiable object agrees with the standard one of
  exponentiable map, in the sense of \cite[p.~7]{joyalmoerdijk95}.
\end{rema}

Given an exponentiable object $E$ in \fibre{C}{U} and an object $A$ in
\fibre{C}{I}, the \emph{canonical cocone from $E$ to $A$} is in the
internal language the cocone of those morphisms from $E$ to
$A$. Formally, it is described as the internal diagram
$(\inct{K}^A,L^A)$, where the internal category $\inct{K}^A$ and the
diagram object $L^A$ are defined as follows. $K^A_0$ is the object
$A^E$, with arrows $s$ and $t$ as in \refdiag{exponential}, and
$K^A_1$ is the pullback
\diag{
  K^A_1 \ar[r]^-{d_0} \ar[d]_x & K^A_0 \ar[d]^s \\
  E^E \ar[r]_-{\overline{t}} & U,
}
where \diag{U & \ar[l]_-{\overline{s}} E^E \ar[r]^-{\overline{t}} &U}
is the exponential of $E$ with itself. In the fibres over $A^E$ and
$E^E$ we have generic maps \func{\varepsilon}{s^*E}{t^*A} and
\func{\overline{\varepsilon}}{\overline{s}^*E}{\overline{t}^*E},
respectively.
  
The codomain map $d_0$ of $\inct{K}^A$ is the top row of the pullback
above, whereas $d_1$ is induced by the composite
\[(\overline{s}x)^*E \xrightarrow{x^*\overline{\varepsilon}}
  (\overline{t}x)^*E \cong (sd_0)^*E
  \xrightarrow{d_0^*\varepsilon} (td_0)^*A\]
via the universal property of $A^E$ and $\varepsilon$.
  
The internal diagram $L^A$ is now the object $s^*E$ in
\fibre{C}{K^A_0}, and the arrow from $d_1^*L^A$ to $d_0^*L^A$ is
(modulo the coherence isomorphisms) $x^*\overline{\varepsilon}$.

When the colimit of the canonical cocone from $E$ to $A$ is $A$
itself, we can think of $A$ as being generated by the maps from $E$ to
it. Therefore, it is natural to introduce the following terminology.

\begin{defi}{generatingobject}
  The object $E$ is called a \emph{generating object} if, for
  any $A$ in $\ct{C}=\fibre{C}{1}$, $A=\colim_{\inct{K}^A}L^A$.
\end{defi}

Later, we shall see how $F$-coalgebras form an indexed category. Then,
a generating object for this category will provide, in the internal
language, a ``generating family'' of coalgebras. The set-theoretic
argument then goes on by taking the coproduct of all coalgebras in
that family. This provides a weakly terminal coalgebra. Categorically,
the argument translates to the following result.

\begin{prop}{weaklyterm}
  Let \indct{C} be an \base{S}-cocomplete \base{S}-indexed category
  with a generating object $E$ in \fibre{C}{U}. Then,
  $\ct{C}=\fibre{C}{1}$ has a weakly terminal object.
\end{prop}

\begin{proof}
  We build a weakly terminal object in \ct{C} by taking the internal
  colimit $Q$ of the diagram $(\inct{K},L)$ in \indct{C}, where
  $K_0=U$, $K_1=E^E$ (with domain and codomain maps $\overline{s}$ and
  $\overline{t}$, respectively), $L=E$ and the map from $d_0^*L$ to
  $d_1^*L$ is precisely $\overline{\varepsilon}$.

  Given an object $A=\colim_{\inct{K^A}}{L^A}$ in \ct{C}, notice that
  the serially commuting diagram
  \diag{
    K^A_1 \ar@<-.8ex>[r]_-{d_0} \ar@<.8ex>[r]^-{d_1} \ar[d]_x &
      K^A_0 \ar[d]^s \\
    E^E \ar@<-.8ex>[r]_-{\overline{t}} \ar@<.8ex>[r]^-{\overline{s}} & U
  }
  defines an internal functor \func{J}{\inct{K}^A}{\inct{K}}.
  We have a commuting triangle of internal \base{S}-categories
  \diag{
    \inct{K}^A \ar[rr]^-J \ar[rd] && \inct{K} \ar[ld]\\
    & 1.
  }
  Taking left adjoint along the reindexing functors which this induces
  on categories of internal diagrams, we get that
  $\colim_{\inct{K}^A}\cong \colim_{\inct{K}}\circ\Lan{J}{}$. Hence, to
  give a map from $A=\colim_{\inct{K}^A}L^A$ to $Q=\colim_{\inct{K}}L$
  it is sufficient to give a morphism of internal diagrams from
  $(\inct{K},\Lan{J}{L^A})$ to $(\inct{K},L)$, or, equivalently, from
  $(\inct{K}^A,L^A)$ to $(\inct{K}^A,J^*L)$, but the reader can easily
  check that these two diagrams are in fact the same.
\end{proof}

Once the coproduct of coalgebras in the ``generating family'' is
formed, the set-theoretic argument is concluded by quotienting it by
its largest bisimulation. One way to build such a bisimulation
constructively is to identify a generating family of bisimulations and
then taking their coproduct.

This suggests that we apply \refprop{weaklyterm} twice; first in the
indexed category of coalgebras, in order to obtain a weakly terminal
coalgebra $(G,\gamma)$, and then in the (indexed) category of
bisimulations over $(G,\gamma)$. To this end, we need to prove
cocompleteness and existence of a generating object for these
categories. The language of inserters allows us to do that in a
uniform way.

Instead of giving the general definition of an inserter in a
2-category, we describe it here explicitly for the 2-category of
\base{S}-indexed categories.

\begin{defi}{inserter}
  Given two \base{S}-indexed categories \indct{C} and \indct{D} and
  two parallel \base{S}-indexed functors
  \func{F,G}{\indct{C}}{\indct{D}}, the \emph{inserter
  $\indct{I}=\indins{F}{G}$ of $F$ and $G$} has as fibre \fibre{I}{X}
  the category whose objects are pairs $(A,\alpha)$ consisting of an
  object $A$ in \fibre{C}{X} and an arrow in \fibre{D}{X} from $F^XA$
  to $G^XA$, an arrow \func{\phi}{(A,\alpha)}{(B,\beta)} being a map
  \func{\phi}{A}{B} in \fibre{C}{X} such that $G^X(\phi)\alpha=\beta
  F^X(\phi)$.
  
  The reindexing functor for a map \func{f}{Y}{X} in \base{S} takes an
  object $(A,\alpha)$ in \fibre{I}{X} to the object
  $(f^*A,f^*\alpha)$,  where $f^*\alpha$ has to be read modulo the
  coherence isomorphisms of \indct{D}, but we shall ignore these
  thoroughly.

  There is an indexed \emph{forgetful} functor
  \func{U}{\indins{F}{G}}{\indct{C}} which takes a pair $(A,\alpha)$
  to its carrier $A$; the maps $\alpha$ determine an indexed natural
  transformation $FU\rTo GU$. The triple $(\indins{F}{G},U,FU\rTo GU)$
  has a universal property, like any good categorical construction,
  but we will not use it in this paper. The situation is depicted as
  below:
  \diaglab{inserter}{\indins{F}{G}\ar[r]^-{U}&\indct{C}\ar@<.8ex>[r]^-F
  \ar@<-.8ex>[r]_-G & \indct{D}.}
\end{defi}

A tedious but otherwise straightforward computation, yields the proof
of the following:

\begin{lemm}{inscocompl}
  Given an inserter as in \refdiag{inserter}, if \indct{C} and
  \indct{D} are \base{S}-cocomplete and $F$ preserves indexed
  colimits, then \indins{F}{G} is \base{S}-cocomplete and $U$
  preserves colimits (in other words, $U$ creates colimits). In
  particular, \indins{F}{G} has all internal colimits, and $U$
  preserves them.
\end{lemm}

\begin{exam}{inserters}
  We shall be interested in two particular inserters, during our
  work. One is the indexed category \indcoal{F}\ of coalgebras for an
  indexed endofunctor $F$ on \indct{C}, which can be presented as the
  inserter
  \diaglab{coalgebras}{
    \indins{\Id}{F}\ar[r]^-{U}&\indct{C}\ar@<.8ex>[r]^-{\Id}
    \ar@<-.8ex>[r]_-F & \indct{C}.}%
  More concretely, $(\indcoal{F})^I=\coalg{F^I}$ consists of pairs
  $(A,\alpha)$ where $A$ is an object and \func{\alpha}{A}{F^IA} a map
  in \fibre{C}{I}, and morphisms from such an $(A,\alpha)$ to a pair
  $(B,\beta)$ are morphisms \func{\phi}{A}{B} in \fibre{C}{I} such
  that $F^I(\phi)\alpha=\beta\phi$. The reindexing functors are the
  obvious ones.
  
  The other inserter we shall need is the indexed category
  \indspan{M}{N} of spans over two objects $M$ and $N$ in \fibre{C}{1}
  of an indexed category. This is the inserter
  \diaglab{span}{
    \indins{\Delta}{<M,N>}\ar[r]^-{U}&\indct{C}\ar@<.8ex>[r]^-{\Delta}
    \ar@<-.8ex>[r]_-{<M,N>} & \indct{C}\prod\indct{C}}
  Where $\indct{C}\prod\indct{C}$ is the product of \indct{C} with
  itself (which is defined fibrewise), $\Delta$ is the diagonal
  functor (also defined fibrewise), and $<M,N>$ is the pairing of the
  two constant indexed functors determined by $M$ and $N$. By this we
  mean that an object in \ct{C} is mapped to the pair $(M,N)$ and an
  object in \fibre{C}{X} is mapped to the pair $(X^*M,X^*N)$.
\end{exam}

\begin{rema}{examcocompl}
  Notice that, in both cases, the forgetful functors preserve
  \base{S}-indexed colimits in \indct{C}, hence both \indcoal{F} and
  \indspan{M}{N} are \base{S}-cocomplete, and also internally
  cocomplete, if \indct{C} is.
\end{rema}

In order to apply \refprop{weaklyterm} to our indexed categories, we
will need to find a generating object for them. This will be achieved
by means of the following two lemmas.

First of all, consider an \base{S}-indexed inserter
$\indct{I}=\indins{F}{G}$ as in \refdiag{inserter}, such that $F$
preserves exponentiable objects. Then, given an exponentiable object $E$ in
\fibre{C}{U}, we can define an arrow $\overline{U}\xrightarrow{r}U$ in
\base{S}\ and an object $(\overline{E},\overline{\varepsilon})$ in
\fibre{I}{\overline{U}}, as follows.

We form the generic map \func{\varepsilon}{s^*F^UE}{t^*G^UE} associated
to the exponential of $F^UE$ and $G^UE$ (which exists because $F$
preserves exponentiable objects), and then define $\overline{U}$ as
the equaliser of the following diagram
\diaglab{mapr}{\overline{U}\ar[r]^-e
  &(G^U\!E)^{F^U\!E} \ar@<.8ex>[r]^-s \ar@<-.8ex>[r]_-t & U,}
the arrow \func{r}{\overline{U}}{U} being one of the two equal
composites $se=te$.

We then put $\overline{E}=r^*E$ and
\diag{
  \overline{\varepsilon}\ =\ F^{\overline{U}}(r^*E) \ar[r]^-{\cong}
  & e^*s^*F^UE \ar[r]^-{e^*\varepsilon}
  & e^*t^*G^UE \ar[r]^-{\cong} & G^{\overline{U}}(r^*E).}
The pair $(\overline{E},\overline{\varepsilon})$ defines an object in
\fibre{I}{\overline{U}}.

\begin{lemm}{existexp}
  The object $(\overline{E},\overline{\varepsilon})$ is exponentiable
  in \indins{F}{G}.
\end{lemm}

\begin{proof}
  Consider an object $(A,\alpha)$ in a fibre \fibre{I}{X}. Then, we
  define the exponential
  $(A,\alpha)^{(\overline{E},\overline{\varepsilon})}$ as follows.

  First, we build the exponential
  \diag{\overline{U} & \ar[l]_-s A^{\overline{E}} \ar[r]^-t & X}
  of $A$ and $\overline{E}$ in \indct{C}, with generic map
  \func{\chi}{s^*\overline{E}}{t^*A}. Because $F$ preserves exponentiable
  objects, we can also form the exponential in \indct{D}
  \diag{\overline{U} & \ar[l]_-{\overline{s}}
        G^X\!A^{F^{\overline{U}}\overline{E}}
        \ar[r]^-{\overline{t}} & X}
  with generic map \func{\overline{\chi}}
  {\overline{s}^*F^{\overline{U}}\overline{E}}{\overline{t}^*G^X\!A}.  By the
  universal property of $\overline{\chi}$, the two composites in
  \fibre{D}{A^{\overline{E}}}
  \[\xymatrix@C+1ex{s^*F^{\overline{U}}\overline{E} \ar[r]^-{\cong} &
    F^{A^{\overline{E}}}s^*\overline{E} \ar[r]^-{F^{A^{\overline{E}}}\chi} &
    F^{A^{\overline{E}}}(t^*A) \ar [r]^-{\cong} &
    t^*F^XA \ar[r]^-{t^*\alpha} & t^*G^XA}\]
  and
  \[\xymatrix@C+1ex{s^*F^{
    \overline{U}}\overline{E} \ar[r]^-{s^*\overline{\varepsilon}} &
    s^*G^{\overline{U}}\overline{E} \ar[r]^-{\cong} &
    G^{A^{\overline{E}}}s^*\overline{E} \ar[r]^-{G^{A^{\overline{E}}}\chi} &
    G^{A^{\overline{E}}}t^*A \ar[r]^-{\cong} & t^*G^XA}\]
  give rise to two maps
  \func{p_1,p_2}{A^{\overline{E}}}{G^XA^{F^{\overline{U}}\overline{E}}}
  in \base{S}, whose equaliser $i$ has as domain the exponential
  $(A,\alpha)^{(\overline{E},\overline{\varepsilon})}$.

  The generic map
  $(si)^*(\overline{E},\overline{\varepsilon})\rTo(ti)^*(A,\alpha)$ in
  \fibre{I}{(A,\alpha)^{(\overline{E},\overline{\varepsilon})}}
  associated to this exponential forms the central square of the
  following diagram, and this commutes because its outer sides are the
  reindexing along the maps $p_1i=p_2i$ of the generic map
  $\overline{\chi}$ above:
  \diag{
    (si)^*F^{\overline{U}}\overline{E}
      \ar[rr]^-{(si)^*\overline{\varepsilon}} \ar[d]_{\cong} &&
    (si)^*G^{\overline{U}}\overline{E} \ar[d]^{\cong} \\
    F^{(A,\alpha)^{(\overline{E},\overline{\varepsilon})}}(si)^*\overline{E}
      \ar[rr]^-{(si)^*(\overline{E},\overline{\varepsilon})}
      \ar[d]_{F^{(A,\alpha)^{(\overline{E},\overline{\varepsilon})}}i^*\chi} &&
    G^{(A,\alpha)^{(\overline{E},\overline{\varepsilon})}}(si)^*\overline{E}
      \ar[d]^{G^{(A,\alpha)^{(\overline{E},\overline{\varepsilon})}}i^*\chi} \\
    F^{(A,\alpha)^{(\overline{E},\overline{\varepsilon})}}(ti)^*A
      \ar[rr]^-{(ti)^*(A,\alpha)} \ar[d]_{\cong} &&
    G^{(A,\alpha)^{(\overline{E},\overline{\varepsilon})}}(ti)^*A
      \ar[d]^{\cong} \\
    (ti)^*F^XA \ar[rr]_-{(ti)^*\alpha} && (ti)^*G^XA.}
  The verification of its universal property is a lengthy but
  straightforward exercise.
\end{proof}

Next, we find a criterion for the exponentiable object
$(\overline{E},\overline{\varepsilon})$ to be generating.

\begin{lemm}{expisgen}
  Consider an inserter of \base{S}-indexed categories as in
  \refdiag{inserter}, where \indct{C} and \indct{D} are
  \base{S}-cocomplete, and $F$ preserves \base{S}-indexed colimits. If
  $(\overline{E},\overline{\varepsilon})$ is an exponentiable object
  in \fibre{I}{\overline{U}} and for any $(A,\alpha)$ in \fibre{I}{1}
  the equation
  \[\colim_{\inct{K}^{(A,\alpha)}}UL^{(A,\alpha)}\cong U(A,\alpha)=A\]
  holds, where $(\inct{K}^{(A,\alpha)},L^{(A,\alpha)})$ is the
  canonical cocone from $(\overline{E},\overline{\varepsilon})$ to
  $(A,\alpha)$, then $(\overline{E},\overline{\varepsilon})$ is
  generating in \indins{F}{G}.
\end{lemm}

\begin{proof}
  Recall from \reflemm{inscocompl} that \indins{F}{G} is internally
  cocomplete and the forgetful functor
  \func{U}{\indins{F}{G}}{\indct{C}} preserves internal
  colimits. Therefore, given an arbitrary object $(A,\alpha)$ in
  \fibre{I}{1}, we can always form the colimit
  $(B,\beta)=\colim_{\inct{K}^{(A,\alpha)}}L^{(A,\alpha)}$. All we need
  to show is that $(B,\beta)\cong(A,\alpha)$. The isomorphism between
  $B$ and $A$ exists because, by the assumption,
  \[B=U(B,\beta)=U\colim_{\inct{K}^{(A,\alpha)}}L^{(A,\alpha)}\cong
    \colim_{\inct{K}^{(A,\alpha)}}UL^{(A,\alpha)}\cong A.\]
  Now, it is not too hard to show that the transpose of the composite
  \[\colim_{\inct{K}^{(A,\alpha)}}F^{\overline{U}}U^{\overline{U}}
    L^{(A,\alpha)} \cong FU\colim_{\inct{K}^{(A,\alpha)}}L^{(A,\alpha)}
    \xrightarrow{\beta}GU\colim_{\inct{K}^{(A,\alpha)}}L^{(A,\alpha)}\]
  is (modulo isomorphisms preserved through the adjunction
  $\colim_{\inct{K}^{(A,\alpha)}}\ladj{\indct{K}^{(A,\alpha)}}^*$) the
  transpose of $\alpha$. Hence, $\beta\cong\alpha$ and we are done.
\end{proof}

As an example, we can show the following result about the indexed
category of spans:

\begin{prop}{genspan}
  Given an \base{S}-cocomplete indexed category \indct{C} and two
  objects $M$ and $N$ in \fibre{C}{1}, if \indct{C} has a generating
  object, then so does the indexed category of spans
  $\indct{P}=\indspan{M}{N}$.
\end{prop}

\begin{proof}
  Recall from \refexam{inserters} that the functor
  \func{U}{\indspan{M}{N}}{\indct{C}} creates indexed and internal
  colimits. If $E$ in \fibre{C}{U} is a generating object for
  \indct{C}, then, by \reflemm{existexp} we can build an exponentiable
  object
  \diag{(\overline{E},\overline{\varepsilon})=M & 
    \ar[l]_-{\overline{\varepsilon}_1} \overline{E}
    \ar[r]^-{\overline{\varepsilon}_2} & N}
  in \fibre{P}{\overline{U}}. We are now going to prove that
  \indspan{M}{N} meets the requirements of \reflemm{expisgen} to show
  that $\overline{E}$ is a generating object.

  To this end, consider a span \diag{(A,\alpha) =
  M&\ar[l]_-{\alpha_1}A\ar[r]^-{\alpha_2}&N} in \fibre{P}{1}. Then, we
  can form the canonical cocone
  $(\inct{K}^{(A,\alpha)},L^{(A,\alpha)})$ from
  $(\overline{E},\overline{\varepsilon})$ to $(A,\alpha)$ in
  \indspan{M}{N}, and the canonical cocone $(\inct{K}^A,L^A)$ from $E$
  to $A$ in \indct{C}. The map \func{r}{\overline{U}}{U} of
  \refdiag{mapr} induces an internal functor
  \func{u}{\inct{K}^{(A,\alpha)}}{\inct{K}^A}, which is an
  isomorphism. Therefore, the induced reindexing functor
  \func{u^*}{\indct{C}^{\inct{K}^A}}{\indct{C}^{\inct{K}^{(A,\alpha)}}}
  between the categories of internal diagrams in \indct{C} is also an
  isomorphism, and hence
  $\colim_{\inct{K}^{(A,\alpha)}}u^*\cong\colim_{\inct{K}^A}$. Moreover,
  it is easily checked that $u^*L^A=UL^{(A,\alpha)}$. Therefore, we have
  $$\colim_{\inct{K}^{(A,\alpha)}}UL^{(A,\alpha)}\cong
    \colim_{\inct{K}^{(A,\alpha)}}u^*L^A\cong
    \colim_{\inct{K}^A}L^A\cong A$$
  and this finishes the proof.
\end{proof}

\section{Final coalgebra theorems}
\label{sec:finalcoalgebra}

In this section, we are going to use the machinery of
\refsec{inserters} in order to prove an indexed final coalgebra
theorem. We then introduce the notion of a class of small maps for a
Heyting pretopos with an (indexed) natural number object, and apply
the theorem in order to derive existence of final coalgebras for
various functors in this context. In more detail, we shall show that
every small map has an M-type, and that the functor \spower\ has a
final coalgebra.

\subsection{An indexed final coalgebra theorem}

In this section, \ct{C} is a category with finite limits and stable
finite colimits (that is, its canonical indexing \indct{C} is a
\base{C}-cocomplete \base{C}-indexed category), and $F$ is an indexed
endofunctor over it (we shall write $F$ for $F^1$). Recall from
\refrema{examcocompl} that the indexed category \indcoal{F}\ is
\base{C}-cocomplete (and the indexed forgetful functor $U$ preserves
indexed colimits).

We say that $F$ is \emph{small-based} whenever there is an
exponentiable object $(E,\varepsilon)$ in \coalg{F^U} such that, for
any other $F$-coalgebra $(A,\alpha)$, the canonical cocone
$(\inct{K}^{(A,\alpha)},L^{(A,\alpha)})$ from $(E,\varepsilon)$ to
$(A,\alpha)$ has the property that
\begin{equation}\label{eq:smallbased}
  \colim_{\inct{K}^{(A,\alpha)}}UL^{(A,\alpha)}\cong U(A,\alpha)=A.
\end{equation}

It is immediate from \refexam{inserters} and \reflemm{expisgen} that,
whenever there is a pair $(E,\varepsilon)$ making $F$ small-based,
this is automatically a generating object in \indcoal{F}. We shall
make an implicit use of this generating object in the proof of:

\begin{theo}{finalcoalg}
  Let $F$ be a small-based indexed endofunctor on a category \ct{C} as
  above. If $F^1$ takes pullbacks to weak pullbacks, then $F$ has an
  indexed final coalgebra.
\end{theo}

Before giving a proof, we need to introduce a little technical lemma:

\begin{lemm}{coalcoeq}
  If $F=F^1$ turns pullbacks into weak pullbacks, then every pair of
  arrows
  \diag{(A,\alpha)\ar[r]^-{\phi}&(C,\gamma)&\ar[l]_-{\psi}(B,\beta)}
  can be completed to a commutative square by the arrows
  \diag{(A,\alpha)&\ar[l]_-{\mu}(P,\chi)\ar[r]^-{\nu}&(B,\beta)}
  in such a way that the underlying square in \ct{C} is a
  pullback. Moreover, if $\psi$ is a coequaliser in \ct{C}, then so is
  $\mu$.
\end{lemm}

\begin{proof}
  We build $P$ as the pullback of $\psi$ and $\phi$ in
  $\ct{C}=\fibre{C}{1}$. Then, since $F$ turns pullbacks into weak
  pullbacks, there is a map \func{\chi}{P}{FP}, making both $\mu$ and
  $\nu$ into coalgebra morphisms. The second statement follows at once
  by the assumption that finite colimits in \ct{C} are stable.
\end{proof}

{\setlength{\parindent}{0pt}{\bf Proof of \reftheo{finalcoalg}.}
  Because \indcoal{F}\ is \base{C}-cocomplete, it is enough, by
  \reflemm{indterm}, to show that the fibre over $1$ of this indexed
  category admits a terminal object.

  Given that $(E,\varepsilon)$ is a generating object in \indcoal{F},
  \refprop{weaklyterm} implies the existence of a weakly terminal
  $F$-coalgebra $(G,\gamma)$. The classical argument now goes on
  taking the quotient of $(G,\gamma)$ by the maximal bisimulation on
  it, in order to obtain a terminal coalgebra. We do that as
  follows. Let $\indct{B}=\indspan{(G,\gamma)}{(G,\gamma)}$ be the
  indexed category of spans over $(G,\gamma)$, i.e.\
  bisimulations. Then, by \refrema{examcocompl}, \indct{B}\ is a
  \base{C}-cocomplete \base{C}-indexed category, and by
  \refprop{genspan} it has a generating object. Applying again
  \refprop{weaklyterm}, we get a weakly terminal span (i.e.\ a weakly
  terminal bisimulation)
  \diag{
    (G,\gamma)&\ar[l]_-{\lambda}(B,\beta)\ar[r]^-{\rho}&(G,\gamma).}
  We now want to prove that the coequaliser
  \diag{
    (B,\beta)\ar@<.8ex>[r]^-{\lambda}\ar@<-.8ex>[r]_-{\rho}&
    (G,\gamma)\ar[r]^-q & (T,\tau)}
  is a terminal $F$-coalgebra.

  It is obvious that $(T,\tau)$ is weakly terminal, since $(G,\gamma)$
  is. On the other hand, suppose $(A,\alpha)$ is an $F$-coalgebra and
  \func{f,g}{(A,\alpha)}{(T,\tau)} are two coalgebra morphisms; then,
  by \reflemm{coalcoeq}, the pullback $s$ (resp.\ $t$) in \ct{C} of
  $q$ along $f$ (resp.\ $g$) is a coequaliser in \ct{C}, which carries
  the structure of a coalgebra morphism into $(A,\alpha)$. One further
  application of \reflemm{coalcoeq} to $s$ and $t$ yields a
  commutative square in \coalg{F}
  \diag{(P,\pi) \ar[r]^-{s'} \ar[d]_{t'} & \bullet \ar[d]^t \\
        \bullet \ar[r]_-s & (A,\alpha)}
  whose underlying square in \ct{C} is a pullback. Furthermore, the
  composite $d=ts'=st'$ is a regular epi in \ct{C}, hence an
  epimorphism in \coalg{F}.

  Write $\widetilde{s}$ (resp.\ $\widetilde{t}$) for the composite of
  $t'$ (resp.\ $s'$) with the projection of the pullback of $f$
  (resp.\ $g$) and $q$ to $G$. Then, the triple
  $((P,\pi),\widetilde{s},\widetilde{t})$ is a span over $(G,\gamma)$;
  hence, there is a morphism of spans
  \[\func{\chi}{((P,\pi),\widetilde{s},\widetilde{t})}
    {((B,\beta),\lambda,\rho).}\]
  It is now easy to compute that $fd=q\lambda\chi=q\rho\chi=gd$, hence
  $f=g$, and the proof is complete.
\mbox{} \hfill $\Box$ \mbox{}\\}

As a particular instance of \reftheo{finalcoalg}, we can recover
the classical result from Aczel \cite[p.~87]{aczel88}.

\begin{corollary}[Final Coalgebra Theorem]
  Any standard functor (on the category of classes) that preserves
  weak pullbacks has a final coalgebra.
\end{corollary}

\begin{proof}
  First of all, notice that preservation of weak pullbacks is
  equivalent to our requirement that pullbacks are mapped to weak
  pullbacks. Moreover, the category of classes has finite limits and
  stable finite colimits. As an exponentiable object, we take the
  class $V$ of all small sets.

  Now, consider a standard functor $F$ on classes (in Aczel's
  terminology). This can easily be seen as an indexed endofunctor,
  since for any two classes $X$ and $I$, one has $X/I\cong X^I$ (so,
  the action of $F$ can be defined componentwise). It is now
  sufficient to observe that every $F$-coalgebra is the union of its
  small subcoalgebras, therefore the functor is small-based in our
  sense.
\end{proof}

\begin{remark}\rm
  With a bit of effort, the reader can see in the present proof of
  \reftheo{finalcoalg} an abstract categorical reformulation of the
  classical argument given by Aczel in his book \cite{aczel88}. In
  order for that to work, he had to assume that the functor preserves
  weak pullbacks (and so did we, in our reformulation). Later, in a
  joint paper with Nax Mendler \cite{aczelmendler89}, they gave a
  different construction of final coalgebras, which allowed them to
  drop this assumption. A translation of that argument in our setting,
  would reveal that the construction relies heavily on the exactness
  properties of the ambient category of classses. Since the functors
  in our examples always preserve weak pullbacks, we preferred
  sticking to the original version of the result (thus making weaker
  assumptions on the category \ct{C}), without boring the reader with
  a (presently unnecessary) second version, which, however, we believe
  can be proved.

  More recently, the work of Ad\'amek, et al.\
  \cite{adamekmiliusvelebil04} has shown that every endofunctor on the
  category of classes is small-based, thereby proving that it has a
  final coalgebra (by Aczel and Mendler's result). Their proof makes a
  heavy use of set theoretic machinery, which would be interesting to
  analyse in our setting.
\end{remark}

\subsection{Small maps}

We are now going to consider on \ct{C} a class of {\em small
maps}. This will allow us to show that certain polynomial functors, as
well as the powerclass functor, are small-based, and therefore we
will be able to apply \reftheo{finalcoalg} to obtain a final coalgebra
for them.

From now on, \ct{C} will denote a Heyting pretopos with an {\em
(indexed) natural number object}. That is, an object \NN, together
with maps \func{0}{1}{\NN} and \func{s}{\NN}{\NN} such that, for any
object $P$ and any pair of arrows \func{f}{P}{Y} and \func{t}{P\prod
Y}{Y}, there is a unique arrow \func{\overline{f}}{P\prod\NN}{Y} such
that the following commutes:
\diag{
  P\prod 1 \ar[d]_{\cong} \ar[r]^-{\id\times 0} &
    P\prod\NN \ar[d]|{<p_1,\overline{f}>} \ar[r]^-{\id\times s} &
    P\prod\NN \ar[d]^{\overline{f}}\\
  P \ar[r]_-{<\id,f>} & P\prod Y \ar[r]_-t & Y.
}
It then follows that each slice $\ct{C}/X$ has a natural number object
$X\prod\NN\rTo X$ in the usual sense. Notice that such categories
have all finite colimits, and these are stable under pullback.

There are various axiomatisations for a class of small maps, starting
with that of Joyal and Moerdijk \cite{joyalmoerdijk95}. In this paper,
we follow the formulation of Awodey et al.\
\cite{awodeybutzsimpsonstreicher04}. A comparison between the two will
appear in \refrema{comparison} below.  A class \smallmap{S}\ of arrows
in \ct{C} is called a class of {\em small maps} if it satisfies the
following axioms:
\begin{description}
\item [(S1)] \smallmap{S}\ is closed under composition and identities;
\item [(S2)] if in a pullback square
             \diag{A \ar[r] \ar[d]_g & B \ar[d]^f\\
                   C \ar[r] & D}
             $f\in\smallmap{S}$, then $g\in\smallmap{S}$;
\item [(S3)] for every object $C$ in \ct{C}, the diagonal
             \func{\Delta_C}{C}{C\prod C} is in \smallmap{S};
\item [(S4)] given an epi \func{e}{C}{D} and a commutative triangle
             \diag{C\ar[rr]^-e\ar@/_/[rd]_f&&D\ar@/^/[ld]^g\\&A,}
             if $f$ is in \smallmap{S}, then so is $g$;
\item [(S5)] if \func{f}{C}{A} and \func{g}{D}{A} are in \smallmap{S},
             then so is their copairing
             \[\func{[f,g]}{C+D}{A}.\]
\end{description}
We say that an arrow in \smallmap{S}\ is {\em small}. We call $X$ a
{\em small object} if the unique map $X\rTo 1$ is small. A {\em small
subobject} $R$ of an object $A$ is a subobject $\xymatrix@1{R\ar@{{
>}->} [r]&A}$ in which $R$ is small. A {\em small relation} between
objects $A$ and $B$ is a subobject $\xymatrix@1{R\ar@{{ >}->}
[r]&A\prod B}$ such that its composite with the projection on $A$ is
small (notice that this does not mean that $R$ is a small subobject of
$A\prod B$).

On a class of small maps, we also require representability of small
relations by means of a {\em powerclass} object:
\begin{description}
\item [(P1)] for any object $C$ in \ct{C}\ there is an object
            $\spower(C)$ and a natural correspondence between maps
            $I\rTo\spower(C)$ and small relations between $I$ and $C$.
\end{description}
In particular, the identity on $\spower(C)$ determines a small
relation $\in_C\subseteq\spower(C)\prod C$. We think of $\spower(C)$ as
the object of all small subobjects of $C$; the relation $\in_C$ then
becomes the membership relation between elements of $C$ and small
subobjects of $C$. The association $C\mapsto\spower(C)$ defines a
covariant functor (in fact, a monad) on \ct{C}. We further require the
two following axioms:
\begin{description}
\item [(I)] The natural number object \NN\ is small;
\item [(R)] There exists a {\em universal small map} \func{\pi}{E}{U}
            in \ct{C}, such that any other small map \func{f}{A}{B}
            fits in a diagram
            \diag{A\ar[d]_f&\bullet\ar[l]\ar[d]\ar[r]&E\ar[d]^{\pi}\\
                  B&\bullet\ar[r]\ar[l]^-{q}&U}
            where both squares are pullbacks and $q$ is epi.
\end{description}
It can now be proved that a class \smallmap{S} satisfying these axioms
induces a class of small maps on each slice $\ct{C}/C$. Moreover, the
reindexing functor along a small map \func{f}{C}{D} has a right
adjoint \func{\Pi_f}{\ct{C}/C}{\ct{C}/D}. In particular, it follows
that all small maps are exponentiable in \ct{C}. 

\begin{rema}{comparison}
  The axioms that we have chosen for our class of small maps subsume
  all of the Joyal-Moerdijk axioms in
  \cite[pp.~6--8]{joyalmoerdijk95}, except for the collection axiom
  $\mathbf{(A7)}$. In particular, the Descent Axiom
  $\mathbf{(A3)}$ can be seen to follow from axioms
  $\mathbf{(S1)-\mathbf(S5)}$ and $\mathbf{(P)}$.

  Conversely, the axioms of Joyal and Moerdijk imply all of our axioms
  except for $\mathbf{(S3)}$ and $\mathbf{(I)}$. Our results in
  \refsec{cst} will imply that, by adding these axioms, a model of the
  weak set theory {\bf CZF$_0$} can be obtained in the setting of
  \cite{joyalmoerdijk95}.
\end{rema}

\subsection{Final coalgebras in categories with small maps}
\label{sec:fincoals}

From now on, we shall consider on \ct{C} a class of small maps
\smallmap{S}. Using their properties, we are now going to prove the
existence of the M-type for every small map \func{f}{D}{C}, as
well as the existence of a final $\spower$-coalgebra.

Let us recall from \cite{demarchivandenberg04} that an exponentiable
map \func{f}{D}{C} in a cartesian category \ct{C} induces on it a {\em
polynomial endofunctor} $P_f$, defined by
\[P_f(X)=\sum_{c\in C}X^{D_c}.\]
Its final coalgebra, when it exists, is called the {\em M-type}
associated to $f$. In fact, the functor $P_f$ is the component over
$1$ of an {\em indexed polynomial endofunctor}, still denoted by
$P_f$, which can be presented as the composite $P_f=\Sigma_C\Pi_fD^*$
of three indexed functors. By this presentation, it follows at once
that $P_f$ preserves pullbacks. The {\em indexed M-type} of $f$ is by
definition the indexed final coalgebra of $P_f$.

In the proof of the following theorem, we will make heavy use of the
internal language of \ct{C}. There, we see $f$ as a signature,
consisting of one term constructor for any $c\in C$ of arity $D_c$,
the fibre of $f$ over $c$. A $P_f$-coalgebra consists of an object $X$
together with a map \func{\gamma}{X}{P_f(X)}, which takes $x\in X$ to
a pair $(c,t)$, where $c\in C$ and $t$ goes from $D_c$ to $X$. The
final $P_f$-coalgebra will then represent the object of all trees
(both well-founded and non-well-founded) over the signature defined by
$f$.

\begin{theo}{finalPfcoalg}
  If \func{f}{D}{C} is a small map in \ct{C}, then $f$ has an
  (indexed) M-type.
\end{theo}
\begin{proof}
  In order to obtain an (indexed) final $P_f$-coalgebra, we want to
  apply \reftheo{finalcoalg}, and for this, what remains to be checked
  is that $P_f$ is small-based. To this end, we first need to find an
  exponentiable coalgebra $(\overline{E},\overline{\varepsilon})$, and
  then to verify condition \refeq{smallbased}.

  The universal small map \func{\pi}{E}{U} in \ct{C} is exponentiable,
  as we noticed after the presentation of axiom {\bf (R)}. Hence,
  unwinding the construction preceding \reflemm{existexp}, we obtain
  an exponentiable object in \indcoal{P_f}. Using the internal
  language of \ct{C}, we can describe
  $(\overline{E},\overline{\varepsilon})$ as follows.

  The object $\overline{U}$ on which $\overline{E}$ lives is described
  as
  \[\overline{U}=\{(u\in U,\func{t}{E_u}{P_f(E_u)})\},\]
  where $E_u$ is the fibre of $\pi$ over $u\in U$, and $\overline{E}$
  is now defined as
  \[\overline{E}=\{(u\in U,\func{t}{E_u}{P_f(E_u)},e\in E_u)\}.\]
  The coalgebra structure \func{\overline{\varepsilon}}{\overline{E}}
  {P_f^{\overline{U}}\overline{E}} takes a triple $(u,t,e)$ (with
  $te=(c,r)$) to the pair $(c,\func{s}{D_c}{\overline{E}})$, where the
  map $s$ takes an element $d\in D_c$ to the triple $(u,t,r(d))$.

  Given a coalgebra $(A,\alpha)$, the canonical cocone from
  $(\overline{E},\overline{\varepsilon})$ to it takes the following
  form. The internal category $\inct{K}^{(A,\alpha)}$ is given by
  \begin{eqnarray*}
    K_0^{(A,\alpha)} & = & 
      \{(u\in U,t:E_u\to P_f(E_u),m:E_u\to A)
      \mid P_f(m)t=\alpha m\};\\
    K_1^{(A,\alpha)} & = & 
      \{(u,t,m,u',t',m',\phi:E_u\to E_{u'}) \mid
      (u,t,m),(u',t',m')\in K_0^{(A,\alpha)}, \\
      && \quad t'\phi=P_f(\phi)t\, \textrm{ and }\, m'\phi=m\}.
  \end{eqnarray*}
  (Notice that, in writing the formulas above, we have used the
  functor $P_f$ in the internal language of \ct{C}; we can safely do
  that because the functor is indexed. We shall implicitly follow the
  same reasoning in the proof of \reftheo{cstafa} below, in order to
  build an (indexed) final \spower-coalgebra.)

  The diagram $L^{(A,\alpha)}$ is specified by a coalgebra over
  $K_0^{(A,\alpha)}$, but for our purposes we only need to consider
  its carrier, which is
  \[UL^{(A,\alpha)}=\{(u,t,m,e)\mid (u,t,m)\in K_0^{(A,\alpha)}
    \textrm{ and }e\in E_u\}.\]
  Condition \refeq{smallbased} says that the colimit of this internal
  diagram in \ct{C} is $A$, but this is implied by the conjunction of
  the two following statements, which we are now going to prove:
  \begin{enumerate}
  \item \label{con1}
        For all $a\in A$ there exists $(u,t,m,e)\in UL^{(A,\alpha)}$
        such that $me=a$;
  \item \label{con2}
        If $(u_0,t_0,m_0,e_0)$ and $(u_1,t_1,m_1,e_1)$ are elements of
        $UL^{(A,\alpha)}$ such that $m_0e_0=m_1e_1$, then there exist
        $(u,t,m,e)\in UL^{(A,\alpha)}$ and coalgebra morphisms
        \func{\phi_i}{E_u}{E_{u_i}} $(i=0,1)$ such that $m_i\phi_i=m$
        and $\phi_i e=e_i$.
  \end{enumerate}
  Condition $\ref{con2})$ is trivial: given $(u_0,t_0,m_0,e_0)$ and
  $(u_1,t_1,m_1,e_1)$, \reflemm{coalcoeq} allows us to fill a square 
  \diag{(P,\gamma)\ar[r]\ar[d] & (E_{u_0},t_0)\ar[d]^{m_0}\\
        (E_{u_1},t_1) \ar[r]_-{m_1} & (A,\alpha),}
  in such a way that the underlying square in \ct{C} is a pullback
  (hence, $P$ is a small object). Therefore, $(P,\gamma)$ is
  isomorphic to a coalgebra $(E_u,t)$, and, under this isomorphism,
  the span
  \diag{(E_{u_0},t_0)&\ar[l](P,\gamma)\ar[r]&(E_{u_1},t_1)}
  takes the form
  \diag{(E_{u_0},t_0)&\ar[l]_-{\phi_0}(E_u,t)\ar[r]^-{\phi_1}&
        (E_{u_1},t_1).}
  Moreover, since $m_0e_0=m_1e_1$, there is an $e\in E_u$ such that
  $\phi_ie=e_i$. Then, defining $m$ as any of the two composites
  $m_i\phi_i$, the element $(u,t,m,e)$ in $UL^{(A,\alpha)}$ satisfies
  the desired conditions.

  As for condition $\ref{con1})$, fix an element $a\in A$. We build a
  subobject $<a>$ of $A$ inductively, as follows:
  \begin{eqnarray*}
    <a>_0 & = & \{a\};\\
    <a>_{n+1} & = & \bigcup_{a'\in<a>_n}t(D_c)\,\textrm{ where }\,
                    \alpha a'=(c,\func{t}{D_c}{A}).\\
  \end{eqnarray*}
  Then, each $<a>_n$ is a small object, because it is a small-indexed
  union of small objects. For the same reason (since, by axiom {\bf
  (I)}, \NN\ is a small object) their union
  $<a>=\bigcup_{n\in\NN}<a>_n$ is small, and it is a subobject of
  $A$. It is not hard to see that the coalgebra structure $\alpha$
  induces a coalgebra $\alpha'$ on $<a>$ (in fact, $<a>$ is the
  smallest subcoalgebra of $(A,\alpha)$ containing $a$, i.e.\ the
  subcoalgebra {\em generated} by $a$), and, up to isomorphism, this
  is a coalgebra \func{t}{E_u}{P_fE_u}, with embedding
  \func{m}{E_u}{A}.  Via the isomorphism $E_u\cong<a>$, the element
  $a$ becomes an element $e\in E_u$ such that $me=a$. Hence, we get
  the desired 4-tuple $(u,t,m,e)$ in $UL^{(A,\alpha)}$.

  This concludes the proof of the theorem.
\end{proof}

\begin{theo}{finalPscoalg}
  The powerclass functor \spower\ has an (indexed) final
  coalgebra.
\end{theo}
\begin{proof}
  It is easy to check that \spower\ is the component on $1$ of an
  indexed functor, and that it maps pullbacks to weak pullbacks.

  Therefore, once again, we just need to verify that \spower\ is
  small-based. We proceed exactly like in the proof of
  \reftheo{finalPfcoalg} above, except for the construction of the
  coalgebra $(<a>,\alpha')$ generated by an element $a\in A$ in
  $\ref{con1})$. Given a \spower-coalgebra $(A,\alpha)$, we construct
  the subcoalgebra of $(A,\alpha)$ generated by $a$ as follows. First,
  we define inductively the subobjects
  \begin{eqnarray*}
    <a>_0 & = & \{a\};\\
    <a>_{n+1} & = & \bigcup_{a'\in<a>_n}\alpha(a').\\
  \end{eqnarray*}
  Each $<a>_n$ is a small object, and so is their union
  $<a>=\bigcup_{n\in\NN}<a>_n$. The coalgebra structure $\alpha'$ is
  again induced by restriction of $\alpha$ on $<a>$.
\end{proof}

\section{The final $\spower$-coalgebra as a model of AFA}
\label{sec:cst}

Our standing assumption in this section is that $\ct{C}$ is a Heyting
pretopos with an (indexed) natural number object and a class
\smallmap{S}\ of small maps. In the last section, we proved that in this
case the $\spower$-functor has a final coalgebra in $\ct{C}$. Now, we
will explain how this final coalgebra can be used to model various set
theories with the Anti-Foundation Axiom. First we work out the case
for the weak constructive theory {\bf CZF$_0$}, and then we
indicate how the same method can be applied to obtain models for
stronger, better known or classical set theories.

Our presentation of {\bf CZF$_0$} follows that of Aczel and Rathjen
in \cite{aczelrathjen01}; the same theory appears under the name of
{\bf BCST*} in the work of Awodey and Warren in
\cite{awodeywarren05}. It is a first-order theory whose underlying logic
is intuitionistic; its non-logical symbols are a binary relation
symbol $\epsilon$ and a constant $\omega$, to be thought of as
membership and the set of (von Neumann) natural numbers,
respectively. Two more symbols will be added for sake of readability,
as we proceed to state the axioms. In order to mark the distinction
between the membership relation of the set theory and that induced by
the powerclass functor inside the category, we shall denote the former
by $\epsilon$ and the latter by the already seen $\in$.

The axioms for {\bf CZF$_0$} are (the universal closures) of the
following statements:
\begin{description}
\item[(Extensionality)] $\forall\;\! x\,(x\;\!\epsilon\;\! a\leftrightarrow x\;\!\epsilon\;\! b\:\!)
                         \rightarrow a=b$
\item[(Pairing)] $\exists\;\!t\,(z\;\!\epsilon\;\! t\leftrightarrow(z=x\lor z=y))$
\item[(Union)] $\exists\;\! t\, (z\;\!\epsilon\;\! t\leftrightarrow \exists\;\! y\,
                (z\;\!\epsilon\;\! y\land y\;\!\epsilon\;\! x))$
\item[(Emptyset)] $\exists\;\! x\, (z\;\!\epsilon\;\! x\leftrightarrow\bot)$
\item[(Intersection)] $\exists\;\! t\, (z\;\!\epsilon\;\! t\leftrightarrow
                       (z\;\!\epsilon\;\! a\land z\;\!\epsilon\;\! b\:\!))$
\item[(Replacement)] $\forall\;\! x\;\!\epsilon\;\! a\,\exists!\;\!y\:\phi\rightarrow
                      \exists\;\! z\, (y\;\!\epsilon\;\! z\leftrightarrow\exists\;\!
                      x\;\!\epsilon\;\! a\: \phi)$
\end{description}
Two more axioms will be added, but before we do so, we want to point
out that all instances of $\Delta_0$-separation follow from these
axioms, i.e. we can deduce all instances of
\begin{description}
\item[($\Delta_0$-Separation)] $\exists\;\! t\, (x\;\!\epsilon\;\! t\leftrightarrow
                                (x\;\! \epsilon\;\! a \land \phi))$
\end{description}
where $\phi$ is a formula in which $t$ does not occur and all
quantifiers are bounded. Furthermore, in view of the above axioms, we
can introduce a new constant $\emptyset$ to denote the empty set, and
a function symbol $s$ which maps a set $x$ to its ``successor'' $x
\cup \{ x \}$. This allows us to formulate concisely our last axioms:
\begin{description}
\item[(Infinity-1)] $\emptyset\;\!\epsilon\;\!\omega\land\forall\;\! x\;\!\epsilon\;\!\omega\,
                      (sx\;\! \epsilon\;\!\omega)$
\item[(Infinity-2)] $\psi(\emptyset)\land\forall\;\! x\;\!\epsilon\;\! \omega\,
                     (\psi(x)\rightarrow\psi(sx))
                     \rightarrow \forall \;\!x\;\!\epsilon\;\! \omega\; \psi(x)$.
\end{description}

It is an old observation by Rieger \cite{rieger57} that models for set
theory can be obtained as fixpoints for the powerclass functor. The
same is true in the context of algebraic set theory (see \cite{butz03}
for a similar result).

\begin{theo}{fixPsmst}
  Every \spower-fixpoint in \ct{C} provides a model of {\bf CZF$_0$}.
\end{theo}
\begin{proof}
  Suppose we have a fixpoint \func{E}{V}{\spower V}, with inverse
  $I$. We call $y$ the \emph{name} of a small subobject $A \subseteq
  V$, when $E(y)$ is its corresponding element in $\spower(V)$. We
  interpret the predicate $x \epsilon y$ as an abbreviation of the
  sentence $x\in E(y)$ in the internal language of \ct{C}. Then, the
  validation of the axioms for {\bf CZF$_0$} goes as follows.

  Extensionality holds because two small subobjects $E(x)$ and $E(y)$
  of $V$ are equal if and only if, in the internal language of \ct{C},
  $z\in E(x)\leftrightarrow z\in E(y)$. The pairing of two elements
  $x$ and $y$ represented by two arrows $1\rTo V$, is given by $I(l)$,
  where $l$ is the name of the (small) image of their copairing
  \func{[x,y]}{1+1}{V}. The union of the sets contained in a set $x$
  is interpreted by applying the multiplication of the monad $\spower$
  to $(\spower E)(E(x))$.  The intersection of two elements $x$ and
  $y$ in $V$ is given by $I(E(x)\cap E(y))$, where the intersection is
  taken in $\spower(V)$. The least subobject $0\subseteq V$ is small,
  and its name \func{\emptyset}{1}{V} models the empty set.

  For the Replacement axiom, consider $a$, and suppose that for every
  $x\;\!\epsilon\;\! a$ there exists a unique $y$ such that
  $\phi$. Then, the subobject $\{y\mid\exists\;\!x\;\!\epsilon\;\!
  a\;\phi\}$ of $V$ is covered by $E(a)$, hence small. Applying $I$ to
  its name, we get the image of $\phi$.

  Finally, the Infinity axioms follow from the axiom {\bf (I)}. The
  morphism \func{\emptyset}{1}{V}, together with the map
  \func{s}{V}{V} which takes an element $x$ to $x\cup\{x\}$, yields a
  morphism \func{\alpha}{\NN}{V}. Since \NN\ is small, so is the image
  of $\alpha$, as a subobject of $V$, and applying $I$ to its name we
  get an $\omega$ in $V$ which validates the axioms Infinity-1 and
  Infinity-2.
\end{proof}

The theorem shows how every $\spower$-fixpoint models a very basic set
theory. Now, imposing extra properties on a fixpoint, we can deduce
the validity of further axioms. For example, in \cite{joyalmoerdijk95}
it is shown how the initial $\spower$-algebra (which is a fixpoint,
afterall!) models the Foundation Axiom. Here, we show how the final
$\spower$-coalgebra satisfies the Anti-Foundation Axiom. To formulate
this axiom, we define the following notions. A (directed) graph
consists of a pair of sets $(n, e)$ such that $n
\subseteq e \times e$. A colouring of such a graph is a function $c$
assigning to every node $x \epsilon n$ a set $c(x)$ such that
\[ c(x) = \{ c(y) \, | \, (x,y)\;\! \epsilon\;\! e \}. \]
This can be formulated solely in terms of $\epsilon$ using the
standard encoding of pairs and functions. In ordinary set theory
(with classical logic and the Foundation Axiom), the only graphs that
have a colouring are well-founded trees and these colourings are then
necessarily unique.

The Anti-Foundation Axiom says:
\begin{description}
\item[(AFA)] Every graph has a unique colouring.
\end{description}

\begin{prop}{mstafa}
  If \ct{C} has an (indexed) final \spower-coalgebra, then this is a
  model for the theory {\bf CZF$_0$}{\rm +}{\bf AFA}.
\end{prop}
\begin{proof}
  We clearly have to check just {\bf AFA}, since any final coalgebra
  is a fixpoint. To this end, note first of all that, because $(V,E)$
  is an indexed final coalgebra, we can think of it as a final
  \spower-coalgebra in the internal logic of \ct{C}.

  So, suppose we have a graph $(n,e)$ in $V$. Then, $n$ (internally)
  has the structure of a \spower-coalgebra \func{\nu}{n}{\spower n},
  by sending a node $x\;\!\epsilon\;\! n$ to the (small) set of nodes
  $y\;\!\epsilon\;\! n$ such that $(x,y)\;\!\epsilon\;\! e$. The
  colouring of $n$ is now given by the unique \spower-coalgebra map
  \func{\gamma}{n}{V}.
\end{proof}

By \reftheo{finalPscoalg}, it then follows at once:
\begin{corollary}
  Every Heyting pretopos with a natural number object and class of
  small maps contains a model of {\bf CZF$_0$}{\rm +}{\bf AFA}.
\end{corollary}

This result can be extended to theories stronger than {\bf
CZF$_0$}. For example, to the set theory {\bf CST} introduced by
Myhill in \cite{myhill75}. This theory is closely related to (in fact,
intertranslatable with) {\bf CZF$_0$}+{\bf Exp}, where {\bf Exp} is
(the universal losure of) the following axiom.
\begin{description}
\item[(Exponentiation)] $\exists\;\! t\,(f\;\!\epsilon\;\! t
                         \leftrightarrow\mbox{Fun}(f,x,y))$
\end{description}

Here, the predicate Fun$(f,x,y)$ expresses the fact that $f$ is a
function from $x$ to $y$, and it can be formally written as the
conjunction of $\forall
a\;\!\epsilon\;\!x\,\exists!b\;\!\epsilon\;\!y\;(a,b)\;\!\epsilon\;\!
f$ and $\forall
z\;\!\epsilon\;\!f\,\exists\;\!a\;\!\epsilon\;\!x, b\;\!\epsilon\;\!
y\;(z=(a,b))$.

\begin{theo}{cstafa}
  Assume the class \smallmap{S} of small maps also satisfies
  \begin{description}
  \item[(E)] The functor $\Pi_f$ preserves small maps for any $f$ in
            \smallmap{S}.
  \end{description}
  Then, \ct{C} contains a model of {\bf CST}{\rm +}{\bf AFA}.
\end{theo}
\begin{proof}
  We already saw how the final \spower-coalgebra $(V, E)$ models
  {\bf CZF$_0$}+{\bf AFA}. Now, {\bf (E)} implies that $A^B$ is
  small, if $A$ and $B$ are, so, in $E(y)^{E(x)}$ is always
  small. This gives rise to a small subobject of $V$, by considering
  the image of the morphism that sends a function $f\in E(y)^{E(x)}$
  to the element in $V$ representing its graph. The image under $I$ of
  the name of this small object is the desired exponential $t$.
\end{proof}

Another example of a stronger theory which can be obtained by imposing
further axioms for small maps is provided by {\bf IZF$^-$}, which
is intuitionistic {\bf ZF} without the Foundation Axiom. It is
obtained by adding to {\bf CZF$_0$} the following axioms:
\begin{description}
\item[(Powerset)] $\exists\;\!y\:\forall\;\!z\:(z\;\!\epsilon\;\!y
                   \leftrightarrow \forall\;\! w\;\! \epsilon\;\!
                   z\:(w\;\!  \epsilon\;\! x))$
\item[(Full Separation)] $\exists \;\!y\: \forall\;\! z\: (z\;\!
                          \epsilon\;\! y \leftrightarrow z\;\!
                          \epsilon\;\! x \land \phi)$
\item[(Collection)] $\forall\;\! y\;\! \epsilon\;\! x\: \exists\;\! w\;
                     \phi \rightarrow \exists\;\! z\:\forall\;\!
                     y\;\!\epsilon\;\!  x\:\exists\;\!
                     w\;\!\epsilon\;\! z\;\phi$
\end{description}
(In Full Separation, $y$ is not allowed to occur in $\phi$.)

In the following theorem, we call a commutative square
\diag{A\ar[r]\ar[d]&B\ar[d]^f\\
      C\ar[r]_-g&D}
a {\em quasi-pullback} if the mediating arrow from $A$ to the pullback
of $f$ and $g$ is epic. By now, the proof of the statement should be
routine (if not, the reader should consult \cite{butz03}):

\begin{theo}{izfafa}
  Assume the class of small maps \smallmap{S} also satisfies
  \begin{description}
  \item[(P2)] if $X \rTo B$ belongs to \smallmap{S}, then so does
              $\spower(X \rTo B)$;
  \item[(M)] every monomorphism is small;
  \item[(C)] for any two arrows $p: Y \rTo X$ and $f: X \rTo A$ where
             $p$ is epi and $f$ belongs to \smallmap{S}, there exists a
             quasi-pullback square of the form
             \diag{ Z \ar[d]_g \ar[r] & Y \ar[r]^p & X \ar[d]^f \\
	            B \ar[rr]_h & & A }
	     where $h$ is epi and $g$ belongs to \smallmap{S}.
  \end{description}
  Then, \ct{C} contains a model of {\bf IZF$^-$}{\rm +}{\bf AFA}.
\end{theo}

\begin{corollary}
  If the pretopos \ct{C} is Boolean, then classical logic is also true
  in the model, which will therefore validate {\bf ZF$^-$}{\rm +}{\bf AFA},
  Zermelo-Fraenkel set theory with Anti-Foundation instead of
  Foundation.
\end{corollary}

Finally, we can build a model for a non-well-founded version of
Aczel's set theory {\bf CZF} in the setting of
\cite{moerdijkpalmgren02}.

The set theory {\bf CZF$^-$}{\rm +}{\bf AFA}, studied by M.~Rathjen in
\cite{rathjen03,rathjen04}, is obtained by adding to ${\bf CZF}_0$
the axiom {\bf AFA}, as well as the following:
\begin{description}
\item[(Strong Collection)] $\forall\;\!x\;\!\epsilon\;\!a\:\exists\;\!y\:
     \phi(x,y)\rightarrow\exists\;\!b\:\mbox{B}(x\;\!\epsilon\;\!a,
      y\;\!\epsilon\;\!b)\:\phi(x,y)$
\item[(Subset Collection)] $\exists\;\!c\:\forall\;\!z\:
     (\forall\;\!x\;\!\epsilon\;\!a\:\exists\;\!y\;\!\epsilon\;\!b\:
     \phi(x,y,z)\rightarrow\exists\;\!d\;\!\epsilon\;\!c\:
     \mbox{B}(x\;\!\epsilon\;\!a,y\;\!\epsilon\;\!d)\:\phi(x,y,z)$
\end{description}
Here, $\mbox{B}(x\;\!\epsilon\;\!a,y\;\!\epsilon\;\!b)\:\phi$ abbreviates:
\[ \forall\;\!x\;\!\epsilon\;\!a\:\exists\;\!y\;\!\epsilon\;\!b\:
   \phi\land\forall\;\!y\;\!\epsilon\;\!b\:\exists\;\!x\;\!\epsilon\;\!a\:
   \phi.\] 

In order for a class of small maps to give a model Subset Collection,
the class has to satisfy a rather involved axiom that will be called
{\bf (F)}. In order to formulate it, we need to introduce some
notation. For two morphisms $A \rTo X$ and $B \rTo X$, $M_X(A, B)$
will denote the poset of multi-valued functions from $A$ to $B$ over
$X$, i.e.\ jointly monic spans in $\ct{C}/X$, 
\diag{ A & \ar@{->>}[l] P \ar[r] & B}
with $P \rTo X$ small and the map to $A$ epic. By pullkback, any
\func{f}{Y}{X} determines an order preserving function
\[ \func{f^*}{M_X(A, B)}{M_Y(f^*A, f^*B)}. \]

\begin{theo}{czfafa}
  Assume the class \smallmap{S} of small maps also satisfies {\bf (C)}
  as in \reftheo{izfafa}, and the following axiom:
  \begin{description}
  \item[(F)] for any two small maps $A \rTo X$ and $B \rTo X$, there
             are an epi $p: X' \rTo X$, a small map $f: C \rTo X'$ and
             an element $P \in M_C(f^*p^*A, f^*p^*B)$, such that for
             any $g: D \rTo X'$ and $Q \in M_D(g^*p^*A, g^*p^*B)$,
             there are morphisms $x: E \rTo D$ and $y: E \rTo C$, with
             $x$ epi, such that $x^* Q \geq y^* P$.
  \end{description}
  Then, \ct{C} contains a model of {\bf CZF$^-$}{\rm +}{\bf AFA}.
\end{theo}

\begin{proof}
  Any fixpoint for \spower\ will model Strong Collection in virtue of
  property {\bf (C)} of the class of small maps.

  Because of {\bf (F)}, the fixpoint will also model the Fullness axiom
  of \cite{aczelrathjen01} (where it is proved to be equivalent to
  Subset Collection over {\bf CZF}$_0$ and Strong Collection).
\end{proof}

Up to this point, we have only given recipes for constructing models
of various non-well-founded set theories. The critical reader might
argue that we have not yet exhibited one such model, since we have not
shown any category to satisfy the given axioms.

To conclude the paper, we present several examples of categories that
satisfy our axioms. Of course, this is not the place to study them in
detail, but we would like to give at least a sketchy presentation. For
a more complete treatment, the reader is advised to look at
\cite{joyalmoerdijk95}. A thorough study of the properties of these
models is the subject for future research.

The most obvious example is clearly the category of classes, where the
notion of smallness is precisely that of a class function having as
fibres just sets. This satisfies all the presented axioms.  Along the
same lines, one can consider the category of sets, where the class of
small maps consists of those functions whose fibres have cardinality
at most $\kappa$, for a fixed infinite regular cardinal $\kappa$. This
satisfies axioms {\bf (S1-5)}, {\bf (P1)}, {\bf(I)}, {\bf (R)}, {\bf
(M)} and {\bf (C)}, but not {\bf (E)}. However, if $\kappa$ is also
inaccessible, then {\bf (E)} is satisfied, as well as {\bf (P2)} and
{\bf (F)}.

Consider the topos $\mathrm{Sh}(\ct{C})$ of sheaves over a site
$\mathcal{C}$, with pullbacks and a subcanonical topology. Then, for
an infinite regular cardinal $\kappa$ greater than the number of
arrows in \ct{C}, define the notion of smallness (relative to
$\kappa$) following \cite[Chapter IV.3]{joyalmoerdijk95}. This
satisfies the axioms {\bf (S1-5)}, {\bf (P1)}, {\bf(I)} and {\bf
(R)}. Moreover, if $\kappa$ is inaccessible, it satisfies also {\bf
(P2)}, {\bf (M)}, {\bf (C)}.

Finally, on the effective topos $\mathcal{E}\!f\!f$ \cite{hyland82} one
can define a class of small maps in at least two different ways. For
the first, consider the global section functor
\func{\Gamma}{\mathcal{E}\!f\!f}{\mathcal{S}et}, and fix a regular
cardinal $\kappa$. Then, say that a map \func{f}{X}{Y} is small if it
fits in a quasi-pullback
\diag{P\ar@{->>}[r]\ar[d]_g&X\ar[d]^f\\
  Q\ar@{->>}[r]&Y}
where $P$ and $Q$ are projectives and $\Gamma(g)$ is $\kappa$-small in
$\mathcal{S}et$. With this definition, the class of small maps
satisfies all the basic axioms {\bf (S1-5)}, {\bf (P1)}, {\bf(I)} and
{\bf (R)}, as well as {\bf (C)} and {\bf(M)}. If $\kappa$ is
inaccessible, it also satisfies {\bf (P2)}.

Alternatively, we can define a map to be small if internally its
fibres are quotients of a subobject of the natural number object of
$\mathcal{E}\!f\!f$. This notion of smallness satisfies all the axioms
apart from {\bf (P2)}.

\bibliographystyle{plain}
\bibliography{shortnames}

\end{document}